# A new method for computing number π

Fernando Alonso Zotes, 03/2021

## Abstract

A family of original formulae for computing number π and its proof are presented in this paper. An algorithm is proposed to test the results. The new method for computing π is interesting from a purely academic point of view, but it is presented with no intention of competing with other efficient formulae already in use for decades, such as the Chudnovsky algorithm [Chudnovsky].

## 1 Introduction

The intangible exploration of higher dimensions has been the purpose of many mathematicians for decades. Several equations are well known for computing the hypervolume of a i-dimensional sphere. In this paper, some of these equations are combined together with the purpose of finding a new method for computing number π. This new method is interesting from an academic point of view, however it is not as efficient as some of the other techniques already well known for decades [Beckmann] and worldwide used for computing number π.

## 2 Development of a new method for computing number π

The hypervolume $V_i(R)$ of an hypersphere of radius $R$ in an i-dimensional space is [NIST]:

$$V_i(R) = k_i R^i$$

Equation 1

with $i \in N$ and:

$$k_i = \begin{cases} i\ odd, & \dfrac{\pi^{\frac{i-1}{2}}}{\prod_{j=1/2}^{j=i/2} j} \\ i\ even, & \dfrac{\pi^{i/2}}{(i/2)!} \end{cases}$$

Equation 2

Note the following recursive relationship:

$$k_i = \frac{2\pi}{i} k_{i-2}$$

$$k_1 = 2,\ k_2 = \pi$$

Equation 3

The volume of the i-dimensional sphere can also be computed by splitting into slices and adding up their volume [Math]:

$$V_i(R) = 2 \int_{x=0}^{x=R} V_{i-1}(r(x))\, dx$$

Equation 4

with $r(x) = \sqrt{R^2 - x^2}$. Equation 1 and Equation 4 yield:

$$V_i(R) = 2 \int_{x=0}^{x=R} k_{i-1} r^{i-1}(x)\, dx$$

Equation 5

Equation 5 yields:

$$V_{i-1}(r(x)) = 2 \int_{y=0}^{y=r(x)} k_{i-2} r^{i-2}(y)\, dy$$

Equation 6

with $r(y) = \sqrt{r^2(x) - y^2}$, thus:

$$V_{i-1}(r(x)) = 2 \int_{y=0}^{y=r(x)} k_{i-2} [r(x)^2 - y^2]^{\frac{i-2}{2}} dy$$

Equation 7

Equation 4 and Equation 7 yield:

$$V_i(R) = 4 k_{i-2} \int_{x=0}^{x=R} \left[ \int_{y=0}^{y=r(x)} [r(x)^2 - y^2]^{\frac{i-2}{2}} dy \right] dx$$

Equation 8

According to the generalized binomial theorem of Newton [Coolidge]:

$$(a-b)^k = \sum_{n=0}^{n=\infty} \left[ (-1)^n \frac{\prod_{j=0}^{j=n}(k-j+1)}{n!\,(k+1)} a^{k-n} b^n \right]$$

Equation 9

which converges if $k \in R$, $|b/a| < 1$. Let it be $a = r^2(x)$ and $b = y^2$, so $|y^2/r^2(x)| < 1$ for $x \in (0, R)$ and $y \in (0, r(x))$, thus Equation 9 yields:

$$(r^2(x) - y^2)^{\frac{i-2}{2}}$$
$$= \sum_{n=0}^{n=\infty} \left[ (-1)^n \frac{\prod_{j=0}^{j=n}(\frac{i-2}{2} - j + 1)}{n!\,(\frac{i-2}{2} + 1)} r^{2(\frac{i-2}{2} - n)}(x) y^{2n} \right]$$
$$= r^{i-2}(x) \sum_{n=0}^{n=\infty} \left[ (-1)^n \frac{2 \prod_{j=0}^{j=n}(\frac{i}{2} - j)}{n!\,i} r^{-2n}(x) y^{2n} \right]$$

Equation 10

Integrating Equation 10 yields:

$$\int_{y=0}^{y=r(x)} [r(x)^2 - y^2]^{\frac{i-2}{2}} dy$$
$$= \int_{y=0}^{y=r(x)} r^{i-2}(x) \sum_{n=0}^{n=\infty} \left[ (-1)^n \frac{2 \prod_{j=0}^{j=n}(\frac{i}{2} - j)}{n!\,i} r^{-2n}(x) y^{2n} \right] dy =$$
$$= r^{i-1}(x) \frac{2}{i} \sum_{n=0}^{n=\infty} \left[ (-1)^n \frac{2 \prod_{j=0}^{j=n}(\frac{i}{2} - j)}{n!\,i(2n+1)} \right] = r^{i-1}(x) P_{i-1}$$

Equation 11

with $P_{i-1}$ independent from $x$ and equal to:

$$P_{i-1} = \frac{2}{i} \sum_{n=0}^{n=\infty} \left[ (-1)^n \frac{\prod_{j=0}^{j=n}(\frac{i}{2} - j)}{n!\,(2n+1)} \right]$$

Equation 12

Integrating Equation 11 yields:

$$\int_{x=0}^{x=R} \left[ \int_{y=0}^{y=r(x)} [r(x)^2 - y^2]^{\frac{i-2}{2}} dy \right] dx$$
$$= \int_{x=0}^{x=R} r(x)^{i-1} P_{i-1} dx$$
$$= P_{i-1} \int_{x=0}^{x=R} (R^2 - x^2)^{\frac{i-1}{2}} dx$$

Equation 13

According to Equation 11:

$$\int_{x=0}^{x=R} (R^2 - x^2)^{\frac{i-1}{2}} dx = R^i P_i$$

Equation 14

with:

$$P_i = \frac{2}{i+1} \sum_{n=0}^{n=\infty} \left[ (-1)^n \frac{\prod_{j=0}^{j=n}(\frac{i+1}{2} - j)}{n!\,(2n+1)} \right]$$

Equation 15

Note Equation 12 and Equation 15 are equivalent.

After replacing Equation 13 and Equation 14 in Equation 8 we get:

$$V_i(R) = 4k_{i-2} R^i P_i P_{i-1}$$

Equation 16

Equation 1 and Equation 16 yield:

$$k_i = 4k_{i-2} P_i P_{i-1}$$

Equation 17

Equation 3 and Equation 17 yield:

$$\frac{2\pi}{i} k_{i-2} = 4k_{i-2} P_i P_{i-1}$$

Equation 18

Solving Equation 18 for π:

$$\pi = 2i P_i P_{i-1}$$

Equation 19

Equation 19 represents a family of formulae to compute number π, with $i \in N$ taking any natural value.

## 3 Algorithm

For implementing Equation 19 in a computer and avoiding factorials of high numbers, the expression of $P_i$ and $P_{i-1}$ given by Equation 15 can be written as follows:

$$P_i = 1 + \sum_{n=1}^{n=\infty} \frac{(-1)^n Q_{i,n}}{2n+1}$$

Equation 20

with:

$$Q_{i,n} = \prod_{j=1}^{j=n} \left[ \frac{i+1}{2j} - 1 \right]$$

Equation 21

Note the following recursive relationship:

$$Q_{i,n} = \left( \frac{i+1}{2n} - 1 \right) Q_{i,n-1}$$

$$Q_{i,1} = \frac{i-1}{2}$$

Equation 22

The algorithm takes into account the fact that $Q_{i,n} = 0$ if $n \geq (i+1)/2$ and $i$ is odd. Equation 20 can be written as follows:

$$P_i = \begin{cases} i\ odd, 1 + \sum_{n=1}^{n=\frac{i-1}{2}} \frac{(-1)^n Q_{i,n}}{2n+1} \\ i\ even, 1 + \sum_{n=1}^{n=\infty} \frac{(-1)^n Q_{i,n}}{2n+1} \end{cases}$$

Equation 23

A python script is located at [Alonso], implementing Equation 19, Equation 22 and Equation 23. The code is as follows:

```
class coef_qq:
  def __init__(self,i):
    self.i=i
    self.v=[(i-1)/2.0]

  def compute(self,n):
    if n>len(self.v):
      # qq(n) is saved in position n-1
      self.v.append(self.compute(n-1)*((self.i+1)/2.0/n-1))
    return self.v[n-1]

class coef_pi:
  def __init__(self,i,Ninf):
    self.qq=coef_qq(i)
    # some terms are cancelled for i odd and n>=(i+1)/2, thus limit set to (i-1)/2 if i odd
    if (i % 2) ==0:
      self.limit=Ninf
    else:
      self.limit=min(Ninf,int((i-1)/2))

  def __computeSumTerm(self,n):
    if (n % 2) == 0:
      out=1.0
    else:
      out=-1.0
    out*=self.qq.compute(n)/(2*n+1)
    return out

  def compute(self):
    sum=1
    for n in range(self.limit):
      # sumatory shall run from 1 to Ninf
```

```
        sum+=self.__computeSumTerm(n+1)
      return sum

def compute_Pi(i,N):
 cpi1=coef_pi(i,N)
 cpi2=coef_pi(i-1,N)
 return 2*i*cpi1.compute()*cpi2.compute()

i=17 # any natural number is valid
N=130 # the larger the better
print('%2.15f'%compute_Pi(i,N)+' = computed π')
from numpy import pi; print('%2.15f'%pi+' = real π')
```

In the algorithm, an upper limit *N* for the infinity sum in Equation 23 is defined. Obviously, better results are expected for higher values of *N*.

The impact of *i* and *N* in the algorithm is quite important. For example, for *i = 5* and *N = 3000000* the resulting estimation of number π is accurate for up to 11 decimal digits, whereas for *i = 17* and *N = 130* the accuracy increases to 15, surprisingly. Therefore, with a good selection of *i* a high accuracy can be achieved even at lower values of *N*.

## 4    Conclusions

A family of formulae has been obtained for computing number π. Proof has been presented. An algorithm has been proposed for computing number π based on the presented family of formulae.

## 5    Future work

As shown above, a good choice of the value of *i* may yield accurate results without the need of using high values for *N*. Further investigation of values for *i* and *N* may lead to interesting conclusions regarding the computational efficiency of the proposed algorithm.